\documentstyle[11pt]{article}

\title{A new criterion for finite non-cyclic groups
\thanks{Project supported by the National Natural Science Foundation(Grant No.10171074).}
}
\author{}
\date{}
\begin{document}
\maketitle
\centerline{Wei Zhou}
\par
\centerline{College of Mathematics \& Finance, Southwest China
Normal University , Chongqing 400715} \centerline{People's Republic
of China}
\centerline{E-mail: zh\_great@hotmail.com}
\par
\centerline{and} \centerline{Wujie Shi}
 \centerline{School of
Mathematics, Suzhou University, Suzhou 215006} \centerline{People's
Republic of China}
 \centerline{E-mail:
wjshi@suda.edu.cn} \par \centerline{and} \centerline{Zeyong Duan}
\centerline{College of Mathematics \& Finance, Southwest China
Normal University , Chongqing 400715}\centerline{People's Republic
of China}
\begin{abstract}
 Let $H$ be a subgroup of a group $G$. We say that $H$ satisfies the power condition with
 respect to $G$,
 or $H$ is a power subgroup of $G$, if there exists a non-negative integer $m$ such that
  $H=G^{m}=<g^{m} | g \in G >$. In this note, the following theorem is proved:
Let $G$ be a group and $k$ the number of non-power subgroups of
$G$. Then (1) $k=0$ if and only if $G$ is a cyclic group(theorem
of F. Sz$\acute{a}$sz) ;(2) $0 < k <\infty$ if and only if $G$ is
a finite non-cyclic group; (3) $k=\infty$ if and only if $G$ is a
infinte non-cyclic group. Thus we get a new criterion for the
finite non-cyclic groups.
\end{abstract}
\par
{\bf Keywords:} power subgroup, cyclic group, Dedekind group.
\par
{\bf 2000MR subject classification:  \ \ 20E07, \ \  20E34, \ \ 20D25}
\par
{ \large \bf 1. Introduction}
\par
Let $H$ be a subgroup of group $G$. We say that $H$ satisfies the
power condition with respect to $G$, or $H$ is a power subgroup of
$G$ if there exists a non-negative integer $m$ such that
$H=G^{m}=<g^{m} | \   g \in G >$. On the other hand, if $H \ne
G^{m}$ for all $m$, we say that $H$ is a non-power subgroup of
$G$. If $H$ is a power subgroup of $G$, the power exponent is the
least non-negative integer $m$ such that $H=G^{m}$.
\par
It is obvious that a nontrivial group $G$ has at least two trivial
power subgroups: $\{1\}$ and $G$ itself, and the power exponents
of the two subgroups are $0$ and $1$, respectively.
\par
The power subgroups have some properties. For example, if $H$ is a
power subgroup of group $G$, then $H$ is a full-invariant
subgroup, in particular $H$ is normal in $G$; and if $H$ is a
non-power subgroup of $G$, then a conjugate subgroup $H^{x}(x \in
G)$ is also a non-power subgroup of $G$.
\par
The number of non-trivial power subgroups affects the structure of
the group. In \cite{Szasz}, F.Sz\'{a}sz  proved that $G$ is a
cyclic group if and only if all subgroups of $G$ are power
subgroups. That is, a non-cyclic group contains at least one
non-power subgroup. In this paper, we generalize the above result
of F. Sz$\acute{a}$sz and prove the following theorem:
\par
\textit{{\bf Theorem.} Let $G$ be a group and $k$ the number of
non-power subgroups of $G$. Then}
\par
\textit{ (1) $k=0$ if and only if $G$ is a cyclic group;}
\par
\textit{ (2) $0 < k <\infty$ if and only if $G$ is a finite
non-cyclic group;}
\par
\textit{ (3) $k=\infty$ if and only if $G$ is a infinte non-cyclic
group.}
\par
\textit{{\bf Remark 1.} The conclusion (1) of this theorem is the
theorem of F.Sz\'{a}sz in \cite{Szasz}.}
\par
\textit{{\bf Remark 2.} Note that for some $k$ in the theorem, the
group $G$ may not exist. In \cite{zhou},\cite{zhou2} the author
proved that the case of $k=1$ or $k=2$ do not occur. But there are
groups having exactly 3 non-power subgroups: For example the
quaternion group $Q_{8}$ of order 8, and $Z_{2} \times Z_{2}$ are
groups possessing just 3 non-power subgroups.}
\par
From the theorem the following corollary gives a new criterion
for finite non-cyclic groups.
\par
\textit{{\bf Corollary.} Suppose $G$ is a non-cyclic group. Then
$G$ is finite if and only if $G$ contains only finitely many
non-power subgroups. }
\par
\textit{{\bf Problem. }For any integer $k(k \ge 3)$, does there
exist groups possessing just $k$ non-power subgroups?}
\par
In the proof of the Theorem, we use the structure of the Dedekind
group. A group is called Dedekind group if all its subgroups are
normal. From Theorem 5.3.7 in \cite{Robinson}, we know that $G$ is
a Dedekind group if and only if $G$ is abelian or the direct
product of a quaternion group of order 8, an elementary abelian
2-group and an abelian group with all its elements of odd order.
Notation is standard and may be found for instance
in\cite{Robinson}. In particular, we denote $k$ the number of of
non-power subgroups of $G$.
\par
{ \large \bf 2. Some lemmas}
\par
\textit{ {\bf Lemma 1}. Let $A$ be an abelian power subgroup in a
group $G$. Then the set of power subgroups of $G$ contained in $A$
coincides with the set of power subgroups of $A$. }
\par
{\bf Proof:} If $G=A$, then there it is nothing to prove. Let $G
\ne A $ and let $l$ be the power exponent of $A$. Then $l > 1$ and
$A=G^{l}$. Suppose that $G^{m} \le A$ for $m \ge 1$ and let $d$ be
 the greatest common divisor of $l,m$.
Then $d=lr+ms$ for some integers $r$ and $s$, and
$g^{d}=(g^{l})^{r}(g^{m})^{s} \in A$ for every $g \in G$, so $
G^{d} \le A$. It is obvious that $A=G^{l} \le G^{d}$, so $d=l$ and
$l$ divides $m$. Hence $m=ls_0$ for some $s_0 \ge 0$. Note that
$A=G^l=<g^l|g \in G>$ is abelian. We get
$A^{s_0}=(G^l)^{s_0}=<g^l|g \in G>^{s_0}=<g^{ls_0}|g \in G>=G^m$.
On the other hand, if $A^{n}$ is a power subgroup of $A$, then
$A^{n}=G^{ln}$ is a power subgroup of $G$. The lemma is proved.
\par
\textit{ {\bf Lemma 2.} If $A/N$ is a non-power subgroup of a
factor-group $G/N$, then $A$ is a non-power subgroup of $G$. }
\par
{\bf Proof:} Suppose the result is false. That is, $A$ is a power
subgroup of $G$. Let $A=G^{m}=<g_{i}^{m} | g_{i} \in G,
i=1,2,\cdots , > $. Since $g_{1}^{m} \cdots g_{r}^{m}
N=(g_{1}N)^{m} \cdots (g_{r} N)^{m},$ $ A/N=(G /N)^{m}$, a
contradiction. Hence the lemma is proved.
\par
Now we need a lemma about the structure of the cyclic group.
\par
\textit{ {\bf Lemma 3.} If $G$ is an infinite abelian group all of
whose proper quotient groups are finite, then $G$ is infinite
cyclic. }
\par
{\bf Proof:} Suppose $G$ is not infinite cyclic. We choose an
element $a_{0} \in G$ and $a_{0} \ne 1$. $<a_{0}>$ is a proper
subgroup of $G$ and $<a_{0}> \ne 1$. Then $G/<a_{0}>$ is finite.
So $G/<a_{0}>$ is finitely generated and $G$ is finitely
generated. Then $G$ is a direct product of finitely many cyclic
groups of infinite or prime-power orders. Let $G=<g_{1}> \times
<g_{2}> \cdots \times <g_{n}>$, and $n >1$ by the assumption. If
there exists an $i$, such that $|g_{i}|=\infty$. Then $G/<g_{1}>
\times \cdots \times <g_{i-1}> \times <g_{i+1}> \times \cdots
\times <g_{n}> \cong <g_{i}>$ is finite by the condition. That is
impossible. So $|g_{i}| < \infty ,\forall i$, which makes $G$ be a
finite group, a contradiction. We prove the lemma.
\par
\textit{ {\bf Lemma 4.} If the number of cyclic subgroups of a
group $G$ is finite, then $G$ is finite. }
\par
{\bf Proof:} Firstly, we have $G$ is a torsion group, otherwise
suppose $g \in G$ and $|g|=\infty$, then $<g^n> $ will be
different cyclic subgroups in $G$ with infinitely many different
$n$, a contradiction.
\par
For any element $g$ in $G$, it can generate a cyclic subgroup
$<g>$. Suppose $G$ have $n$ cyclic subgroups and $g_i(i=1,\cdots
,n)$ be the generated elements of all the cyclic subgroups.
$G=\cup ^n_{i=1} <g_i>$. So $|G| \le \sum^n_{i=1} |g_i|<\infty$,
that is to say $G$ is finite.
\par
{\large \bf 3. Proof of the theorem}
\par
{\bf Proof of the theorem:} We need only to prove the necessity of
(2).
\par
Case 1. $G$ is abelian.
\par
i) Suppose $G$ is a torsion abelian group, we prove that $G$ is a
finite group.
\par
Let $\pi(G)$ be the set of all primes dividing the orders of
elements of $G$. For every prime $p$, let $G_{p}$ be the set of
all $p$-elements in $G$. Then $G_{p}$ is a subgroup of $G$, that
is, the $p$-component of $G$. We claim $|\pi(G)|$ is finite. If
$|\pi (G)|=\infty$, then there must exist a prime $p \in \pi(G)$
such that $G_p$ is a power subgroup for there are finitely many
non-power subgroups. Hence there exists a positive integer $m$
such that $G_p=G^m$, and then $G/G_p=G/G^m$, which means the
exponent of $G/G_p$ is finite. But $G/G_p \cong \prod _{q \in
\pi(G) -\{p\}} G_q$ . So the exponent of $G/G_p$ can not be finite
since $|\pi(G)|$ is infinite and then we get the contradiction.
\par
If $G_{p}$ is infinite for some $p$ , then it contains infinitely
many finite cyclic subgroups $C_{i}$ and hence we may assume that
$C_{i}$ is a power subgroup in $G$ and so $G/C_{i}$ is of finite
exponent, hence $G$ is of finite exponent and $G$ is a direct
product of cyclic subgroups of finite order. In particular,
$G_{p}$ is of finite exponent and since $G_{p} \simeq G/N$ for
some subgroup $N$, $G_{p}$ has only finitely many non-power
subgroups in $G_{p}$ by Lemma 2. Since $G_{p}$ has only finitely
many power subgroups, i.e., $G_{p},G_{p}^{p}, \cdots,
G_{p}^{p^{l}}$, where $p^{l}$ is the power exponent of $G_{p}$,
join with its finitely many non-power subgroups, $G_{p}$ has only
finitely many subgroups and so $G_{p}$ is finite. Then for every
$p \in \pi (G)$, $G_{p}$ is finite. So $G$ is finite.
\par
ii) Suppose $G$ is a non-torsion abelian group. Then $G$ is an
infinite abelian group with some elements of infinite order. We
prove$G$ is cyclic, which is contrary to $k > 0$($k$ is the number
of non-power subgroups of $G$).
\par
Firstly, we prove that the subgroup $T$ of $G$ consisting of all
elements of finite orders in $G$ is finite. Evidently, no subgroup
of $T$ can be a power subgroup of $G$, otherwise $G$ would be
torsion and so, by i), $G$ is finite, a contradiction. Hence,
since has only finitely many non-power subgroups, we have that $T$
contains only finitely many subgroups. Thus $T$ have only finitely
many cyclic subgroups and $T$ is finite by Lemma 4.
\par
Now $T \ne G$ and $\overline{G}=G/T$ is an infinite torsion-free
group. By Lemma 2 , $\overline{G}$ has only finitely many
non-power subgroups. We prove $\overline{G}$ is cyclic.
\par
We claim that $\overline{G}/\overline{H}$ is finite for all
non-trivial subgroup $\overline{H}$ of $\overline{G}$. In fact,
let $\overline{H}$ be a nontrivial subgroup of $\overline{G}$. If
$\overline{H}$ is a power subgroup of $\overline{G}$, then there
exists a non-negative integer $m$ such that
$\overline{H}=\overline{G}^{m}$. Now $\overline{G} /\overline{H}$
is torsion, and it has only finitely many non-power subgroups by
Lemma 2. Then $\overline{G}/\overline{H}$ is finite by i). On the
other hand, suppose $\overline{H}$ is a non-power subgroup. Since
$\overline{G}$ is torsion-free with finitely many non-power
subgroups, there must exist $1 \ne \overline{H}_{1}
<\overline{H}$, with $\overline{H}_{1}$ being a power subgroup of
$\overline{G}$. Then $\overline{G}/\overline{H}_{1}$ is finite and
then $\overline{G} /\overline{H}$ is finite.
\par
By Lemma 3, $\overline{G}$ is cyclic. Hence  $G=T \times <z>$
where $z$ is of infinite order. It is obvious that $G^{m}=T^{m}
\times <z^{m}>$, so if $s=|T| >1$ then $T \times <z^{sn}>$ is a
non-power subgroup of $G$ for each $n=2,3,\cdots$. Thus $T=1$ and
so $G$ is cyclic, a contradiction.
\par
Case 2. $G$ is non-abelian.
\par
If $G$ contains no non-power subgroups, then $G$ is cyclic and
every subgroup of a cyclic group is a power subgroup. So, suppose
that $G$ contains some non-power subgroups, and let $H_{1},\cdots
, H_{s}$ be all of those. Since every conjugate of $H_{i}
(i=1,\cdots ,s)$ is also a non-power subgroup of $G$, so $H_{i}$
has only finitely many conjugate subgroups. Then the normalizer
$N_{i}$ of $H_{i}$ has a finite index in $G$ and hence the
subgroup $K_{i}= \cap _{g \in G} g^{-1} N_{i}g$ is normal in $G$
and has a finite index in $G$. Let $K=\cap_{i=1}^{s} K_{i}$, then
$K$ is a normal subgroup of finite index in $G$ and so it
normalizes every $H_{i}$. Observe that every subgroup of $K$ is
normal in $K$. In fact, if $L \le K$ and $L$ is a power subgroup
of $G$, then $L \unlhd G$ and hence $L \unlhd K$; if $L$ is a
non-power subgroup of $G$, then $L=H_{i}$ for some $i=1,\cdots,s$,
and by construction of $K$, $H_{i} \unlhd K$. Thus $K$ is a
Dedekind group. So $K$ is either abelian or a direct product of an
abelian group and the quaternion group of order 8. Hence the
center $Z=Z(K)$ of $K$ is of finite index in $K$ and evidently $Z
\lhd G$. Thus $Z$ is an abelian normal subgroup of finite index,
say $m$, in $G$. In particular, $x^{m} \in Z$ for every $x \in G$.
Let $l$ be the smallest natural number such that $x^{l} \in Z$ for
every $x \in G$. Then $G^{l} \le Z$.
\par
Observe that $G^{l}$ is a subgroup of finite index in $Z$.
Otherwise, then $Z/G^{l}$ is infinite and contains infinitely many
proper non-trivial subgroups $R_{i}/G^{l}, i=1,2,\cdots$. By
assumption, there exists $t$ such that $R_{t}$ is a power subgroup
of $G$, so $R_{t}=G^{r}$ for some $r \ne l$. It is easy to get
$r<l$, but this contradicts the choice of $l$.
\par
Thus the index of $G^{l}$ in $Z$ is finite, so we may replace $Z$
by $G^{l}$ and assume that $Z$ is an abelian power subgroup of
finite index in $G$, $Z=G^{l}$.
\par
Observe that, by Lemma 1, every power subgroup of $G$, being
contained in $Z$, is a power subgroup of $Z$, and every power
subgroup of $Z$ is a power subgroup of $G$. In particular, $Z$
contains only finitely many non-power subgroups, and hence $Z$ is
finite or an infinite cyclic group. If $Z$ is finite, then $G$ is
finite and so the result is proved. Thus we suppose $Z$ is an
infinite cyclic group, that is, $Z=<x>$ with $|x|=\infty$. If
$C_{G}(Z) \ne G$, then there exists $y \in G$ such that
$x^{y}=y^{-1} x y \ne x \in <x>$. Obviously $Z^y=Z$. Since $x$ and
$x^{-1}$ are the only generators of $<x>$, we have $x^{y}=x^{-1}$.
\par
Let $Y_i=<x^{-i}yx^i>, i=1,2,\cdots$. We claim there are infinite
many different $Y_i$. Suppose it is false. Let $Y_{i_1}, \cdots,
Y_{i_s}$ be all the different such type subgroups. Then there are
infinite subgroups of this type equal to $Y_{i_t},1 \le t \le s$.
Without any loss, we suppose
$Y_{i_t}=Y_{m_1}=\cdots=Y_{m_r}=\cdots$. If $|y|=\infty$, by the
generator properties of infinite cyclic group we have
$y^{x^{m_1}}=(y^{x^{m_r}})^{\varepsilon_r}, r=2,3,\cdots$, where
$\varepsilon_r=1$ or $-1$. Hence there must exist $r_0,r_1$ such
that $\varepsilon_{r_0}=\varepsilon_{r_1}$ and $r_0 \ne r_1$. And
then $y^{x^{m_{r_0}}}=y^{x^{m_{r_1}}}$. Suppose $m_{r_0} <m_{r_1}$
 without any loss. Then we have $y=y^{x^{m_{r_1}-m_{r_p}}}$.
But since $x^y=x^{-1}$, we get
$y^{x^{m_{r_1}-m_{r_p}}}=x^{-2(m_{r_1}-m_{r_0})}y$. Therefore
$y=y^{x^{m_{r_1}-m_{r_p}}}=x^{-2(m_{r_1}-m_{r_0})}y$, and then
$x^{-2(m_{r_1}-m_{r_0})}=1$. We get contradiction for
$|x|=\infty$. On the other hand if $|y|=n$. Similarly there are
infinite subgroups such that $Y_{m_1}=\cdots=Y_{m_r}=\cdots$. By
the generator properties of cyclic group we have
$y^{x^{m_1}}=(y^{x^{m_r}})^{\varepsilon_r}, r=2,3,\cdots$, where
$1 \le \varepsilon_r \le n$ is coprime to $n$. Hence there must
exist $r_0,r_1$ such that $\varepsilon_{r_0}=\varepsilon_{r_1}$
and $r_0 \ne r_1$. And then $y^{x^{m_{r_0}}}=y^{x^{m_{r_1}}}$. We
get a contradiction in the same way as before.
\par
So there are infinite many different $Y_i$ which are conjugate
with each other. In particular, these infinite many $Y_{i}$ are
non-normal in $G$ and non-power subgroups of $G$, which is
contrary to the assumption. Thus $C_{G}(Z)=G$ and hence $Z \le
Z(G)$. So the center of $G$ has a finite index in $G$. By Schur
Theorem (\cite{Robinson} Theorem 10.1.3), the commutator subgroup
$C=[G,G]$ is finite. If $C=1$ then $G$ is abelian contrary to the
assumption. So let $y$ be an element of prime order $p$ in $C$.
Obviously $Z \cap <y>=1$ and $<y> \lhd<Z,y>$ by $Z \le Z(G)$, so $
<Z,y>= Z \times <y>$ and all subgroups $Z^{r} \times <y>$ are
distinct for $r=1,2,\cdots $. Thus there exists $r$ such that
$A=Z^{r} \times <y>$ is a power subgroup of $G$. $A$ contains, by
Lemma 1, only a finitely many non-power subgroups of $A$ and hence
$A$ is finite or cyclic. But this is not true. The theorem is
proved.
\par
{\bf Acknowledegement } The authors would like to thank Prof. V.D.
Mazurov for his help.
\par

\end{document}